\documentclass[12pt]{amsart}
\usepackage{amscd}
\def\NZQ{\Bbb}               
\def\NN{{\NZQ N}}

\def\RR{{\NZQ R}}
\def\CC{{\NZQ C}}

\def\frk{\frak}               
\def\aa{{\frk a}}

\def\mm{{\frk m}}

\def\Phi{{\frk n}}
\def\Phi{{\frk N}}
\def\opn#1#2{\def#1{\operatorname{#2}}} 
\opn\chara{char}
\opn\length{\ell}
\opn\pd{pd}
\opn\rk{rk}
\opn\projdim{proj\,dim}
\opn\injdim{inj\,dim}
\opn\rank{rank}
\opn\depth{depth}
\opn\and{and}
\opn\grade{grade}
\opn\height{height}
\opn\embdim{emb\,dim}
\opn\codim{codim}

\opn\Tr{Tr}
\opn\bigrank{big\,rank}
\opn\superheight{superheight}\opn\lcm{lcm}
\opn\trdeg{tr\,deg}%
\opn\reg{reg}
\opn\lreg{lreg}
\opn\ini{in}
\opn\div{div}
\opn\Div{Div}
\opn\cl{cl}
\opn\Cl{Cl}
\opn\Spec{Spec} \opn\Supp{Supp} \opn\supp{supp} \opn\Sing{Sing}
\opn\Ass{Ass} \opn\Min{Min}
\opn\Ann{Ann}
\opn\Rad{Rad}
\opn\Soc{Soc}
\opn\Im{Im}
 \opn\Ker{Ker} \opn\Coker{Coker} \opn\Am{Am}
\opn\Hom{Hom} \opn\Tor{Tor} \opn\Ext{Ext} \opn\End{End}
\opn\Aut{Aut} \opn\id{id}

\opn\nat{nat}
\opn\pff{pf}
\opn\Pf{Pf}
\opn\GL{GL}
\opn\SL{SL}
\opn\mod{mod}
\opn\ord{ord}
\opn\cl{cl}
\opn\conv{conv}
\opn\ext{ext}
\opn\rad{rad}
\opn\red{red}
\opn\aff{aff}
\opn\con{conv}
\opn\relint{relint}
\opn\st{st}
\opn\lk{lk}
\opn\cn{cn}
\opn\core{core}
\opn\vol{vol}
\opn\link{link}
\opn\star{star}
\opn\gr{gr}

\def\pot#1#2{#1[\kern-0.28ex[#2]\kern-0.28ex]}

\opn\dirlim{\underrightarrow{\lim}}
\opn\inivlim{\underleftarrow{\lim}}
\let\sect=\cap

\let\tensor=\otimes
\let\iso=\cong

\let\Sect=\bigcap
\let\Dirsum=\bigoplus

\let\to=\rightarrow

\def\Implies{\ifmmode\Longrightarrow \else
     \unskip${}\Longrightarrow{}$\ignorespaces\fi}
\def\implies{\ifmmode\Rightarrow \else
     \unskip${}\Rightarrow{}$\ignorespaces\fi}
\def\iff{\ifmmode\Longleftrightarrow \else
     \unskip${}\Longleftrightarrow{}$\ignorespaces\fi}

\let\:=\colon
\newtheorem{Theorem}{Theorem}[section]
\newtheorem{Lemma}[Theorem]{Lemma}
\newtheorem{Corollary}[Theorem]{Corollary}
\newtheorem{Proposition}[Theorem]{Proposition}
\newtheorem{Remark}[Theorem]{Remark}
\newtheorem{Remarks}[Theorem]{Remarks}
\newtheorem{Example}[Theorem]{Example}

\let\epsilon\varepsilon
\let\phi=\varphi
\let\kappa=\varkappa
\def\qed{\ifhmode\textqed\fi
   \ifmmode\ifinner\quad\qedsymbol\else\dispqed\fi\fi}
\def\textqed{\unskip\nobreak\penalty50
    \hskip2em\hbox{}\nobreak\hfil\qedsymbol
    \parfillskip=0pt \finalhyphendemerits=0}
\def\dispqed{\rlap{\qquad\qedsymbol}}

\opn\dis{dis}
\def\pnt{{\raise0.5mm\hbox{\large\bf.}}}

\begin{document}

\title{Minimal Monomial Reductions and the Reduced Fiber Ring of an Extremal Ideal}

\author{Pooja Singla}
\address{FB Mathematik, Universit\"at Duisburg-Essen, Campus Essen, 45117 Essen, Germany}
\email{pooja.singla@uni-duisburg-essen.de}

\begin{abstract}
Let $I$ be a monomial ideal in a polynomial ring
$A=K[x_1,\ldots,x_n]$. We call a monomial  ideal $J$ to be a
 minimal monomial reduction ideal of $I$ if there exists no proper
monomial ideal $L \subset J$ such that $L$ is a reduction ideal of
$I$. We prove that there exists a unique minimal monomial
reduction ideal $J$ of $I$ and we show that  the maximum  degree
of a monomial generator of $J$ determines the slope $p$ of the
linear function $\reg(I^t)=pt+c$ for $t\gg 0$. We determine the
structure of the reduced fiber ring $\mathcal{F}(J)_{\red}$ of $J$
and show that $\mathcal{F}(J)_{\red}$ is isomorphic to the inverse
limit of an inverse system  of semigroup rings determined by
convex geometric properties of $J$.
\end{abstract}

\maketitle
\section*{Introduction}

Let $I$ be a monomial ideal in a
polynomial ring $A=K[x_1,\ldots,x_n]$ over a field $K$.
Let $G(I)$ denote the unique minimal monomial set of generators of $I$.

Cutkosky-Herzog-Trung \cite{CHT} and independently Kodiyalam
\cite{KO} have shown that for any graded ideal $I$ in a polynomial
ring $A=K[x_1,\ldots,x_n]$, the regularity of $I^t$ is a linear
function $pt+c$ for large enough $t$. Also the coefficient $p$ of
the linear function  is known and it is given by the
$\min\{\theta(J): J \mbox{ is a graded reduction ideal of}\; I\}$,
see \cite{KO}. Here $\theta(J)$ denotes the maximum of the degrees of
elements  in $G(J)$.

In  Section 2 we give  a convex geometric interpretation for this
coefficient $p$ for any monomial ideal $I \subset A$: let $S$ be
any set of monomials in $A$. We denote by $\Gamma(S)\subset \NN^n$
the set of exponents of the monomials in $S$. Now let $J$ be the
monomial ideal which is determined by the property that
$\Gamma(G(J))= \ext(I)$, where $\ext(I)$ denotes the extreme
points of the convex set $\conv(I)$.  Here $\conv(I)$ denotes the
convex hull of the elements of the set $\Gamma(I)$ in $\RR^n$.
This convex set is commonly called the {\em Newton polygon} of
$I$. We show in Proposition \ref{T1} that the ideal $J$ is the
unique minimal monomial reduction ideal of $I$, that is, there
exists no proper monomial ideal $L \subset J$ such that $L$ is
again a reduction ideal of $I$. It turns out that $p=\theta(J)$.
In other words, $p=\max\{\deg x^a\: a\in\ext(I)\}$.

We call a reduction ideal $L$ of
$I$ to be a {\em Kodiyalam reduction} if $\theta(L)=p$. Thus
the ideal $J$ generated by monomials whose
exponents belong to $\ext(I)$ is a Kodiyalam reduction.

We call a monomial ideal $L$ to be an {\em extremal ideal} if
$\Gamma(G(L))=\ext(L)$. In other words, $L$ is an extremal ideal
if $L$ is its own minimal monomial reduction. Notice  that each
squarefree monomial ideal is an extremal ideal. Let $\mu(L)$
denote the number of generators in a minimal generating set of  a
graded ideal $L$. It is easy to see that $\mu(\Rad I)$ is bounded
above by $|\ext(I)|$ for any monomial ideal $I \subset A$.

In Section 3 we describe  the faces of $\conv(I^m)$ for a monomial
ideal $I$, and compare the supporting hyperplanes and the faces of
$\conv(I^{n_1})$ and $\conv(I^{n_2})$ for two positive integers
$n_1,n_2.$

In Section 4 we determine the structure of the reduced fiber ring
$\mathcal{F}(L)_{\red}$ of an extremal ideal $L$. For any graded
ideal $L \subset A=K[x_1,\ldots,x_n]$, the fiber ring
$\mathcal{F}(L)$ is defined to be
$\mathcal{R}(L)/\mm\mathcal{R}(L)=\Dirsum_{n \geq 0} L^n/mL^n$
where ${\mathcal R}(L)$ is  the Rees ring and
$\mm=(x_1,\ldots,x_n) \subset A$ is the graded maximal  ideal of
$A$. The main motivation to study the structure of the reduced
fiber ring of an extremal ideal is to determine the dimension of
the fiber ring of an arbitrary  monomial ideal. Let $I\subset A$
be a monomial ideal and $J \subset I$ be its minimal monomial
reduction. Then $J$
 is an extremal ideal, and $\dim \mathcal{F}(I)= \dim \mathcal{F}(J)=\dim\mathcal{F}(J)_{\red}$.
So as far as dimension is concerned it is enough to consider the
reduced fiber ring $\mathcal{F}(J)_{\red}$ of the extremal ideal
$J$, whose structure is in general much simpler than that of
$\mathcal{F}(J)$.

Let $\mathcal{F}_c$ denote the set of all compact faces of
$\conv(I)$. It is shown in Lemma \ref{imp} that for  each $F
\in\mathcal{F}_c$, we have $F=\conv\{a_{j_1},\ldots,a_{j_t}\}$
where $F\sect \ext(I)=\{a_{j_1},\ldots,a_{j_t}\}.$ For each $F\in
{\mathcal F}_c$ we put $K[F]=K[x^{a_j}t\:  a_j \in F]$. As the
main result of Section 4 we show in Theorem \ref{limit}  that
${\mathcal F}(J)_{\red} \iso \underleftarrow{\lim}_{F \in
\mathcal{F}_c}K[F]$. As an application to the  Theorem \ref{limit}
we get in the particular case of monomial ideals a result of
Carles Bivia-Ausina \cite{BA} on the analytic spread of a Newton
non-degenerate ideal.


Let $\overline{L}$ denotes the integral closure of an ideal $L$.
In Section 5,  using convex geometric arguments, we  show in
Theorem \ref{ell} that $\overline{I^\ell}=J\overline{I^{\ell-1}}$
where $\ell$ is the analytic spread of $I$.  If we assume that
$I^a$ is integrally closed for $a \leq \ell-1$, then as a
corollary of Theorem \ref{ell}, we obtain that
$I^{\ell}=JI^{\ell-1}$, and that $I$ is a normal ideal.

I am very much grateful to Prof.\ Herzog for many helpful
discussions and comments.

\section{Some Preliminaries on the Convex Geometry of Monomial Ideals}

Let $I$  be a monomial ideal in a polynomial ring
$A=K[x_1,\ldots,x_n]$ over a field $K$. We denote by $G(I)$ the
unique  minimal monomial generating set of $I$.

For a monomial $u=x^a=x_1^{a(1)}\cdots x_n^{a(n)}\in A$ we denote
by $\Gamma(u)$ the exponent vector $(a(1),\ldots, a(n))$ of $u$.
Similarly, if $S$ is any set of monomials in $A$, we set
$\Gamma(S)=\{\Gamma(u)\: u\in S\}$.

We denote the convex hull of $\Gamma(I)$ by $\conv(I)$. Here $\Gamma(I)=\{a: x^a \in I \}$.
Recall
that $\conv(I)$ is a polyhedron. A polyhedron can be defined as
the intersection of finitely many closed half spaces. A polyhedron
may also be thought of as the sum of a polytope (which is the
convex hull of a finite set of points) and the positive cone
generated by a finite set of vectors. Indeed these two notions are
equivalent, (see \cite[Theorem 1.2]{ZG}).

Suppose that $G(I)=\{x^{a_1},\ldots, x^{a_s}\}$, then
$$\conv(I)= \conv\{a_1,a_2,\ldots,a_s\}+ \RR^n_{\geq 0},$$ see \cite[Lemma
4.3]{RV}. Here the positive cone $\RR^n_{\geq 0}$ denotes the set
of vectors $u \in \RR^n$ such that $u(i) \geq 0$ for all $i
=1,\ldots,n$. It follows that $\conv(I)$  is a polyhedron. It is
called the {\em Newton polyhedron} of $I$.

Let $H_i =\{ v \in \RR^n \;|\; \langle v,u_i \rangle = c_i \}$
where $u_i \in \RR^n$ , $c_i \in \RR$ for  $i=1,\ldots,m$  be the
hyperplanes in $\RR^n$ such that $\conv(I)= \{ v \in \RR^n |\;
\langle v,u_i \rangle \geq c_i,\; i=1,\ldots,m\}$. We observe

\begin{Lemma} \label{L1} The vectors $u_i$ belong to $\RR^n_{\geq 0}$ for  $i=1,\ldots,m$.
\end{Lemma}

\begin{proof} We prove that $\langle e_j,u_i \rangle=u_i(j) \geq 0$ for all $i,j$.
Here the vectors $e_j=(0,\ldots,0,1,0,\ldots,0)$
are the canonical unit vectors in $\RR^n$ for $j=1,\ldots,n$ and
$1$ being at $j$th place. Let $a \in \Gamma(I)$, then $a+ te_j \in
\conv(I)$ for all $j$ and $ t\in \RR_{\geq 0}$. Hence
 $ \langle a+ te_j,u_i \rangle \geq c_i$ for all $i,j$. Suppose $\langle e_{j_0}, u_{i_0} \rangle < 0$
 for some $j_0,i_0$.
 Then we have
 $\langle a+te_{j_0}, u_{i_0} \rangle < c_{i_0}$  for $t \gg 0$,  which is a contradiction.
 \end{proof}
Before proceeding further we need to set up some terminology from convex geometry  (see \cite{GB}).
We define the notions of exposed points and extreme points for a convex set $X \subset \RR^n$.
 A point $a \in X $ is said to be an {\em extreme point},
provided all $b,c \in X$, $0< \lambda <1$, and $a= \lambda
b+(1-\lambda)c$ imply $a=b=c.$ We denote this set of extreme
points by $\ext(X)$.

Let $H=\{ v \in \RR^n \;|\; \langle v,u \rangle = c \}$ be a
hyperplane where $ u\in \RR^n, c\in \RR $. We denote by $H_+$ the
nonnegative closed half space defined by $H$, i.e. $H_+=\{ v \in
\RR^n \;|\; \langle v,u \rangle \geq c \}$. We say $H$ is a {\em
supporting hyperplane} of a closed convex set $X,$ if $X \subset
H_+$ and $ X \cap H \neq \emptyset$. Again, we may notice, as in  Lemma \ref{L1} that for every
supporting hyperplane $H  =\{ v \in \RR^n \;|\; \langle v,u
\rangle = c \} \subset \RR^n$ of $\conv(I)$ one has $u \in \RR^n_{
\geq 0}$.

A set $ F  \subset X$ is called a {\em face} of $X$, if either $F=\emptyset $, or $F=X$,
or if there exists a supporting hyperplane
$H$ of $X$ such that $F=X \cap H$. We call $F$ to be a {\em proper face}
of $X$ if $F \neq X$ and $F \neq \emptyset$.

Let $F$ be a proper face of $\conv(I)$. Let
$H=\{ v \in \RR^n \;|\; \langle v,u \rangle = c \}$ be a supporting hyperplane of $\conv(I)$ such that
$F=H \cap \conv(I)$.
It may be observed
that $F$ is a compact  face of $\conv(I)$ if and only if the vector $u \in
(\RR_{+}\backslash\{0\})^n $ i.e. $u(j)>0$ for all $j=1,\ldots,n$.

We now define exposed points of $X$ which we denote by $\exp(X)$.
A point $a \in X$ is called an {\em exposed point} of $X$ if the
set $\{a\}$ consisting of single point is a face of $X$. Hence for
every $a \in \exp(X)$ there exists a supporting hyperplane $H  =\{
v \in \RR^n \;|\; \langle v,u \rangle = c \} \subset \RR^n$ such
that $ \{a\}= X \cap H$ i.e. $\langle a,u \rangle =c$ and $
\langle b,u \rangle > c$ for all $b \in X, b \neq a$.

 We denote the extreme points of $\conv(I)$ by $\ext(I)$
and the exposed points of $\conv(I)$ by $\exp(I)$.

\begin{Proposition} \label{P1} Let $I$  be a monomial ideal in a polynomial
ring $A=K[x_1,\ldots,x_n]$ over a field $K$. Then, $a \in \exp(I)$
implies $x^a \in G(I)$.
\end{Proposition}

\begin{proof} Let $a \in \exp(I)$ and
$H  =\{ v \in \RR^n \;|\; \langle v,u \rangle = c \}$ $\subset \RR^n$ be a supporting hyperplane of $\conv(I)$ such that
 $H \cap \conv(I)=\{a\}$. Notice that $u \in (\RR_+ /\{0\})^n$.

Let $G(I)=\{x^{a_1},\ldots,x^{a_s}\}$. Then  $\conv(I)=\conv \{a_1,a_2,\ldots,a_s\}+ \RR^n_{\geq 0}.$ Therefore
$a=\sum_{i=1}^s k_ia_i + v$ where $\sum_{i=1}^sk_i=1$, $k_i \geq
0$, $v \in \RR^n_{\geq 0}$. Now, since $\langle a_i, u \rangle
\geq c$ and $\langle w,u \rangle > 0$ for any $0 \neq w \in
\RR^n_{\geq 0}$, $\langle a,u \rangle =c$ implies $a=a_i$ for some
$i$.
\end{proof}

\begin{Remark} \label{R1}
{\em For any closed convex set $X \subset \RR^n$, one has $\exp(X)
\subset \ext(X)$ and $\ext(X) \subset \cl(\exp(X))$ where
$\cl(\exp(X))$ denotes the closure of $X$ in $\RR^n$ with respect
to usual topology (see \cite[Statement 3 \mbox{and} 9, Section 2.4
]{GB}). In case $X= \conv(I)$, one has $\exp(I)$ is a finite set.
Therefore $\cl(\exp(I))=\exp(I)$, and  hence $\exp(I)=\ext(I)
\subset \Gamma(G(I))$.}
\end{Remark}

\section{Minimal Monomial Reduction Ideal}

In this section we show that for any monomial ideal $I \in
A=K[x_1,\ldots,x_n]$, there exists a unique minimal monomial
reduction ideal $J$ of $I$. We also show that the minimal monomial
reduction ideal $J$ of a monomial ideal $I$ is a Kodiyalam
reduction of $I$

Let $L \subset A=K[x_1,\ldots,x_n]$ be a graded ideal. An ideal $N
\subset L$ is said to  be a {\em reduction ideal of $L$}, if there
exists a positive integer $m$ such that $NL^{m-1}=L^m$. Let
$\bar{I}$ denote the integral closure of an ideal $I$. It is known
that $N\subset L$ is a reduction ideal of $L$ if and only if
$\overline{N}=\overline{L}$(see \cite[Exercise 10.2.10(c)]{BH}).

 Now let $I \subset K[x_1,\ldots,x_n]$ be a
monomial ideal. We say a monomial ideal $J \subset I$ a {\em
minimal monomial reduction ideal of $I$} if there exists no proper
monomial ideal $J' \subset J$ such that $J'$ is a reduction ideal
of $I$. For a monomial ideal one has
$$\Gamma(\bar{I})=\conv(I)\cap \NN^n$$ (see \cite[Exercise 4.22
]{ED}). Hence a monomial ideal $J\subset I$ is a reduction ideal
of $I$ if and only if $\conv(J)=\conv(I)$.

\begin{Proposition} \label{T1}
Let $I$  be a monomial ideal in a polynomial ring
$A=K[x_1,\ldots,x_n]$ over a field $K$ with
$\ext(I)=\{a_1,\ldots,a_r\}$. Then the  ideal $J=(x^{a_1},\ldots,
x^{a_r})$ is the unique minimal monomial reduction ideal of $I$.
\end{Proposition}

\begin{proof} To show that $J$ is a reduction ideal of $I$ is equivalent to prove that
$\conv(I)=\conv(J)$. For any monomial ideal $L\subset A$, we know
that $\conv(L)=\conv(\Gamma(G(L)))+\RR^n_{\geq 0}$, \cite[Lemma
4.3]{RV}. We also have $\Gamma(G(J))=\ext(I)$. On the other hand
it follows easily from \cite[Section 8.9]{sch} that
\begin{eqnarray}
\label{schade} \conv(I) =\conv(\ext(I))+ \RR^n_{\geq 0}.
\end{eqnarray}
These facts imply that $\conv(I)=\conv(J)$.

Again, it is also easy to see that $J$ is the unique minimal
monomial reduction of $I$. In fact,  let $L$ be any other monomial
reduction ideal of $I$. We show that $J \subset L$. We have
$\conv(I)=\conv(L)$, and so $\ext(L)=\ext(I)$. We know that
$\ext(L)=\exp(L) \subset \Gamma(G(L))$, by Lemma \ref{P1}.
Therefore we have $\Gamma(G(J))\subset \Gamma(G(L))$.
\end{proof}


\noindent
For all nonnegative integers $m$, we define the ideal $J^{[m]} := (x^{ma_1},\ldots,x^{ma_r})$.
\begin{Corollary} \label{last} The ideal $J^{[m]}$ is the
unique minimal monomial reduction ideal of $I^m$ for all $m$.
\end{Corollary}
\begin{proof}  Let us fix an $m$, and denote by $J_m$ the unique monomial reduction ideal of $I^m$.
First notice that $J^{[m]}$ is a monomial reduction ideal of $I^m$.
Indeed, as $J^{[m]}$ is a monomial reduction ideal of $J^m$ and $J^m$ is
a monomial reduction ideal of $I^m$, we have $J^{[m]}$ is a reduction ideal of
$I^m$. Therefore $J_m\subset J^{[m]}$, by Theorem \ref{T1}.

Next we claim that $ \ext(I^m) \supset \{ma_1,\ldots,ma_r\}$, and
this will imply that  $J^{[m]}\subset J_m$, by Theorem \ref{T1}.

Let $ H_i = \{ v \in \RR^n \;|\; \langle v,u_i \rangle = c_i\}$ be
a supporting hyperplane of $\conv(I)$ such that $H_i \cap
\conv(I)=\{a_i\}$ for $i=1,\ldots,r$. We define the hyperplanes
$mH_i=\{ v \in \RR^n \;|\; \langle v,u_i \rangle = mc_i\}$,
$i=1,\ldots,r$ and show that $mH_i$ is a supporting hyperplane of $\conv(I^m)$ with
$mH_i\sect \conv(I^m)=\{ma_i\}$. This
then will imply the above claim.

It is clear that $ma_i\in mH_i\sect \conv(I^m)$. Now let $a\in
\Gamma(I^m)$ be an arbitrary element. Then
$a=\sum_{j=1}^ma_{i_j}+v$ where $v\in \NN^n$.  It follows that
$\langle a, u_i\rangle\geq mc_i$, and is equal to $mc_i$ if and
only if $a=ma_i$, as $u_i\in (\RR_{+}/\{0\})^n$. It follows that
$\langle b,u_i\rangle \geq mc_i$ for all $b\in \conv(I^m)$, and
equality holds if and only if $b=ma_i$.
\end{proof}

Let $I$ be a graded ideal in a polynomial ring
$A=K[x_1,\ldots,x_n]$ over a field $K$. The $i$th regularity of an
ideal $I$ is defined to be $\reg_i(I)=\max\{j :
\Tor_i^A(I,K)_{i+j} \neq 0\}$ and the Castelnuovo-Mumford
regularity of $I$ is defined to be $\reg(I)=\max\{\reg_i(I)-i\}.$

Cutkosky-Herzog-Trung \cite{CHT} and independently Kodiyalam
\cite{KO} have shown that $\reg(I^t)=pt+c$  for $t\gg 0$. Also the
coefficient  of the linear function is known and it is given by
$$p=\min\{\theta(J): J \mbox{ is a graded reduction ideal of}\;
I\},$$ see  \cite{KO}. Here $\theta(J)$ denotes the maximum of the degrees
of elements in $G(J)$. We define a reduction ideal $J$ of $I$
to be a {\em Kodiyalam reduction} if $\theta(J)=p$.

More generally, it is shown in \cite{CHT}  that  $\reg_i(I^t)
=p_it+q_i$ for $t\gg 0$ are linear functions. From the arguments
in Kodiyalam's paper \cite{KO} it follows immediately that
$p_0=p$.

\begin{Corollary} \label{C1} Let $I$ be a monomial ideal in $K[x_1,\ldots,x_n]$, then the minimal monomial
reduction ideal J of $I$
is a
Kodiyalam reduction.
\end{Corollary}
\begin{proof}
The proof is very much on the line of arguments of Kodiyalam (see
\cite[Proposition 4]{KO}). By the very definition of $p$, we have
$\theta(J) \geq p$. We now show that $\theta(J) \leq p$. It is
enough to find a monomial reduction ideal $L$ such that $\theta(L)
\leq p $, as $G(J)\subset G(L)$ because
$\Gamma(G(J))=\ext(I)=\ext(L)\subset \Gamma(G(L))$. Notice that $\ext(I)=\ext(L)$, as $L \subset I$ being
a reduction ideal of  $I$, we have $\conv(I)=\conv(L)$.

Consider the minimal  monomial generating system of $I$, given by $f_1,\ldots,f_s$ where $\deg f_i=d_i$ for all $i$ and
$d_1\leq \cdots \leq d_s$.
Let $j$ be the largest integer such that $f_j^k \notin \mm I^k$ for any $k$ where $\mm$ is the maximal graded ideal in $A$.
 Then $\reg_0(I^t)\geq d_jt$ for all $t$. Set $L=(f_1,\ldots,f_j)$
 and $P=(f_{j+1},\ldots,f_s)$.
Clearly $L$ is a monomial ideal with $\theta(L)=d_j$. We claim that $L$ is a reduction ideal of $I$.
By the very choice of $j$, $ P^t \subset \mm I^t$ for some $t$. Then $I^t=(L+P)^t=L(L+P)^{t-1}+P^t \subset LI^t+\mm I^t$.
Hence by Nakayama's lemma, it follows that $L$ is a reduction ideal of $I$.
Now as $\theta(L)=d_j$ and $d_jt \leq pt+q_0$ for $t \gg 0$. We have $d_j \leq p$. Hence $\theta(L) \leq p$.
\end{proof}
We call a monomial ideal $L$ an {\em extremal ideal},
if $G(L)=\ext(L)$. In other words, a monomial ideal $L$ is an extremal
ideal if it is the minimal monomial reduction of itself. In
particular, the ideal $J$ in Theorem \ref{T1} is an extremal
ideal.

\begin{Remarks}
{\em
1. Every squarefree monomial ideal is an extremal ideal.
Let $N \subset A$ be a squarefree monomial ideal and let $x^a \in G(N)$ be  a monomial generator.
We show that $a \in \ext(N)$.
As $N$ is squarefree, for all $i$, one has $a(i)=1$ or $a(i)=0$.
Let $r \leq n$ be the cardinality of  $i$'s such that $a_i=1$.
We define a vector $u \in \NN^n$ given by $u(i)=1$  if $a(i)=1$ and $u(i)=n+1$ if $a(i)=0$.
We claim that the hyperplane  $S=\{v \in \RR^n : \langle v,u \rangle = r \}$ is a supporting hyperplane of
$\conv(N)$
with $S \cap \conv(N)=\{a\}$, which will imply that $a \in \ext(N)$.
Clearly, $S \cap \conv(N) \supset \{a\}$. Let $b \in \conv(N)=\conv(\Gamma(G(N))+\RR^n_{\geq 0}$ with $b \neq a$
be an arbitrary element. We claim that $\langle b,u \rangle >r$. Notice that it is enough to
consider $b \in \Gamma(G(N))$. Since $x^a,x^b \in G(N)$, we notice that
there exists an $i$ such that $b(i)=1$ and $a(i)=0$. Hence $\langle b,u \rangle \geq n+1$ and so
$\langle b,u \rangle >r$. Hence the claim.

Let $\mu(L)$ denote the number of generators in a minimal
generating set of  a graded ideal $L$.

2. Let $I \subset A$
be a monomial ideal. Then we have $ \mu(\Rad I)\leq |\ext(I)|$.
Infact, let $J \subset I$ be the minimal monomial reduction ideal of $I$.
Then one has $\Rad J=\Rad I$. Hence $\mu(\Rad I)=\mu(\Rad J)\leq \mu(J)=|G(J)|=|\ext(I)|$. }
\end{Remarks}

\section {A description of the faces of $\conv(I^m)$}

Let $I=(x^{a_1},x^{a_2},\ldots,x^{a_s})\subset
A=K[x_1,\ldots,x_n]$ be a monomial ideal. We may assume that
$\ext(I):=\{a_1,\ldots,a_r\}$ is the set of extremal points of the
convex hull of $I$ after a proper rearrangement of generators.
Then $J= ( x^{a_1},x^{a_2},\ldots,x^{a_r}) $ is the minimal
monomial reduction ideal of $I$, see Theorem \ref{T1}.

Next we consider the set of faces of $\conv(I)$.
Let $\mathcal{F}$ denote the set of proper faces and
$\mathcal{F}_c \subset \mathcal{F}$ denote the set of compact
faces of $\conv(I)$. Let $F \in \mathcal{F}$ and $S:=\{ v \in
\RR^n \;|\; \langle v,u \rangle = c\}$ be a supporting hyperplane
of $\conv(I)$ such that $S \cap \conv(I)= F$. It may be observed
that $F \in \mathcal{F}_c$ if and only if the vector $u \in
(\RR_{+}\backslash\{0\})^n $. For $j=1,\ldots,n$, we define
$e_j=(0,\ldots,0,1,\ldots,0)\in\RR^n$ to be the unit vectors,
$1$ being at $j$th place.
\noindent
With this notation, we have

\begin{Lemma} \label{imp}
Let $F\in \mathcal F$ be a face of $\conv(I)$, and let  $S=\{v\in
\RR^n\: \langle v,u\rangle=c\}$ be a  supporting hyperplane of
$\conv(I)$ such that $F=S\sect \con(I)$. Then $F\sect \ext(I)\neq
\emptyset$, and $$F=\conv\{a_{j_1},\ldots,a_{j_t}\}+\sum_{\{j\:
u(j)=0\}}\RR_{\geq 0} e_j,$$ where $F\sect
\ext(I)=\{a_{j_1},\ldots,a_{j_t}\}.$
\end{Lemma}
\begin{proof} Let $a \in \conv(I)$. Then
$a=\sum_i^rk_ia_i+v$ with $\sum k_i=1$, $k_i \geq 0$, $v \in
\RR^n_{\geq 0}$ by  Equation \ref{schade}. Suppose $F \cap \ext(I)= \emptyset$.
Then $\langle a_i, u \rangle>c$  for all
$i=1,\ldots,r$. Therefore we have $\langle a, u \rangle>c$. Hence
$F= S \cap \conv(I)= \emptyset$, a contradiction.

Now let $F\sect \ext(I)=\{a_{j_1},\ldots,a_{j_t}\}.$ First let $F$
be a compact face, then $u \in (\RR_+ \backslash \{0\})^n$. As
$\langle a_i, u \rangle >c$ for all $a_i \in \ext(I)\backslash
\{a_{j_1},\ldots,a_{j_t}\}$ and $\langle v, u \rangle> 0$ for all
$0 \neq v \in \RR^{n}_{\geq 0}$, we notice that $\langle a, u
\rangle=c$ if and only if $a \in \conv\{a_{j_1},\ldots,a_{j_t}\}$.
Hence $F=\conv\{a_{j_1},\ldots,a_{j_t}\}$.

Now let $F$ be an noncompact face and let $Z=\{j : u(j)=0\}$. Notice
that the set $Z\neq \emptyset$. As $\langle a_i, u \rangle >c$ for
all $a_i \in \ext(I)\backslash \{a_{j_1},\ldots,a_{j_t}\}$ and
$\langle v, u \rangle \geq 0$ for all $v \in \RR^{n}_{\geq 0}$
with  $\langle v, u \rangle=0$ if and only if $v \in \sum_{j \in Z
}\RR_{\geq 0} e_j$,  we see that  $\langle a, u \rangle=c$ if and
only if $a \in \conv\{a_{j_1},\ldots,a_{j_t}\}+\sum_{\{j\:
u(j)=0\}}\RR_{\geq 0} e_j$.
\end{proof}

As an immediate consequence of Lemma \ref{imp} we obtain

\begin{Corollary} \label{supporting} Let $S=\{ v \in \RR^n \;|\; \langle v,u \rangle = c\}$ be a hyperplane. Then
$S$ is a supporting hyperplane of $\conv(I)$ if and only if $\langle a_i, u \rangle \geq c$ for all $a_i \in \ext(I)$
and $\langle a_j, u \rangle=c$
for some $a_j \in \ext(I)$.
\end{Corollary}


\begin{Lemma} \label{2.1}
Let  $S =\{ v \in \RR^n \;|\; \langle v,u \rangle = c \}$ where $
u\in \RR^n$, $c\in \RR,$ be a  hyperplane,  and let $n_1, n_2\geq
1$ two integers and  $q=n_2/n_1$. Then $S$ is a supporting hyperplane of $\conv(I^{n_1})$
if and only if $qS=\{ v \in \RR^n \;|\; \langle v,u \rangle = qc \}$
is supporting hyperplane of $\conv(I^{n_2})$.
\end{Lemma}

\begin{proof} We know by Corollary \ref{last} that $\ext(I^m)=(ma_1,\ldots,ma_r)$ for all $m \geq 1$.
Now $S$  is a supporting hyperplane of $\conv(I^{n_1})$ if and
only if  $\langle n_1a_i, u \rangle \geq c$ for all $n_1a_i \in
\ext(I^{n_1})$ and $\langle n_1a_j, u \rangle=c$ for some $n_1a_j
\in \ext(I^{n_1})$. This is the case if and only if $\langle
n_2a_i, u \rangle=\langle (n_2/n_1)n_1a_i, u \rangle=q\langle
n_1a_i, u \rangle  \geq qc$ and $\langle n_2a_j, u
\rangle=q\langle n_1a_j, u \rangle =qc$. This is equivalent to say
that  $qS$ is a supporting hyperplane of $\conv(I^{n_2})$, see
Corollary \ref{supporting}.
\end{proof}

 Let $\mathcal{F}$
be the set of proper faces of $\conv(I)$. For each $F\in
\mathcal{F}$ we choose a hyperplane $S=\{ v \in \RR^n \;|\; \langle v,u \rangle = c\}$   with $F=S\sect
\conv(I)$. Then by Lemma \ref{2.1}, for any nonnegative integer $m$, the hyperplane $mS$ is a supporting
hyperplane of $\conv(I^m)$, and we set $mF=mS\sect \conv(I^m)$. It
is easy to see that this definition does not depend on the choice
of $S$. Indeed,
\[
mF=\conv\{ma_{j_1},\ldots,ma_{j_t}\}+\sum_{\{j\:
u(j)=0\}}\RR_{\geq 0} e_j
\]
if $F\sect \ext(I)=\{a_{j_1},\ldots,a_{j_t}\}$.
 We denote by $m{\mathcal F}$ the set of proper
faces of $\conv(I^m)$.
\noindent
As an immediate consequence of Lemma \ref{2.1} we get

\begin{Corollary}
\label{multiplefaces} The map ${\mathcal F}\to m{\mathcal F}$,
$F\mapsto mF$ is bijective.
\end{Corollary}

\section{The structure of the reduced fiber ring of an extremal
ideal}

The main result of this section is Theorem \ref{limit} which gives
us the structure of the reduced fiber ring of an extremal ideal.
We proceed gradually towards it preparing the ground to prove it.
We will use all the notation from previous section.

Recall that a monomial ideal $L \subset A=K[x_1,\ldots,x_n]$ is
said to be an extremal ideal if $\Gamma(G(L))=\ext(L)$. In other
words an extremal ideal is the minimal monomial reduction of
itself, see Proposition \ref{T1}.

The main motivation to study the structure of the reduced fiber
ring $\mathcal{F}(J)_{\red}$ of an extremal ideal is to determine
the dimension of the fiber ring $\mathcal{F}(I)$ for any monomial
ideal $I$. As one notices that $\dim \mathcal{F}(I)=\dim
\mathcal{F}(J)=\dim \mathcal{F}(J)_{\red}$, therefore it is enough
to consider the reduced fiber ring $\mathcal{F}(J)_{\red}$ as far
as the dimension is concerned. We will see that in general the
structure of the reduced fiber ring of an extremal ideal is more
simple than that of the original fiber ring.

For the proof of Theorem \ref{2.4}  we shall need the following

\begin {Lemma}\label{2.3}
Let $ a=  \sum_{i=1}^r l_ia_i$
 where $l_i$ are nonnegative integers, $\sum l_i=m$ and $\ext(I)=\{a_1,\ldots,a_r\}$.
If $\{ a_i : l_i \neq 0\} \not\subset F$ for some $F \in \mathcal
{F}$, then $a \notin mF$.
\end{Lemma}
\begin{proof} Let $S=\{ v \in \RR^n \;|\; \langle v,u \rangle = c \}$ be a
supporting hyperplane of $\conv(I)$ such that $S \cap \conv(I)=F$.
Then $mS=\{ v \in \RR^n \;|\; \langle v,u \rangle = mc \}$ is  a
supporting hyperplane of $\conv(I^m)$ such that $mS \cap
\conv(I^m)=mF$.

Suppose that $a \in mF$. Then  we have   $\langle a,
 u \rangle=mc$.
 Since $\{ a_i : l_i \neq 0\} \not\subset F$, there exists at least one $ j $
such that $\langle a_j, u \rangle > c$ which implies  $\langle a, u \rangle > mc$, a contradiction.
\end{proof}

\begin{Remark}\label{extra} {\em From the above lemma, it follows that
if the set  $\{ a_i : l_i \neq 0\}
\not\subset F$ for any $F \in \mathcal {F}$, then $a \notin G$ for any $G \in m\mathcal{F}$.
Indeed, as for every $G \in m\mathcal{F}$ there exists $F \in \mathcal{F}$ such that
$G=mF$,
by Corollary \ref{multiplefaces}.}
\end{Remark}


The following theorem is crucial in our study of the structure of
the reduced fiber ring of an extremal ideal.

\begin{Theorem} \label{2.4} Let  $J$ be an extremal ideal with $G(J)=\{f_1,\ldots, f_r\}$ and
$f_j=x^{a_{j}}$ for $j=1,\ldots, r$.
Let $Z=\{a_{j_1},\ldots,a_{j_t}\}$ be a subset of $\Gamma(G(J))$.
 Then the
following conditions are equivalent:
\begin{enumerate}
\item $ Z \subset F$ for some compact face  $F \in \mathcal{F}$;

\item For all $l_i\geq 0$ one has $f_{j_1}^{l_1}\cdots
f_{j_t}^{l_t}\in G(J^m)$  where $m=\sum_{i=1}^tl_{i}$;

\item For all $l_{i}\gg 0$ one has $f_{j_1}^{l_{1}}\cdots
f_{j_t}^{l_{t}}\in G(J^m)$  where $m=\sum_{i=1}^tl_{i}$.
\end{enumerate}
\end{Theorem}

\begin{proof}
$(1) \Rightarrow (2)$ Suppose there exists some nonnegative
integers $l_{i}$ such that $f'=f_{j_1}^{l_1}\cdots f_{j_t}^{l_t}
\notin G(J^{m})$ where $m=\sum l_{i}$. Then there exists  $g \in
G(J^{m})$ such that $f'=hg$ where $\deg h > 0$. Let $S :=\{ v \in
\RR^n \;|\; \langle v,u \rangle = c \}$ be a supporting hyperplane
such that $F=S \cap \conv(J).$ Notice that as $F$ is a compact
face, the vector $u$ belongs to $(\RR_{+}\backslash\{0\})^n.$ Now
since the set $Z \subset F$, $\langle a_{j_k}, u \rangle =c$ for
all $k=1,\ldots,t$.  Then we have $\langle \Gamma(f'), u
\rangle=mc$, but since $\langle \Gamma(h), u \rangle >0$ and
$\langle \Gamma(g), u \rangle \geq mc$,  one has $\langle
\Gamma(hg), u \rangle>mc$, a contradiction.

$(2) \Rightarrow (3)$ is trivial.

$(3) \Rightarrow (1)$ Suppose if $Z \not\subset F$ for any compact
face $F \in \mathcal{F}$, then we prove that for all $l_{i}\gg 0$
we have $f_{j_1}^{l_1}\cdots f_{j_t}^{l_t}\notin G(J^m)$  where
$m=\sum_{i=1}^tl_{i}$.

Let $f=f_{j_1}\cdots f_{j_t}$. We will show that
$f^{m_0}=f_{j_1}^{m_0}\cdots f_{j_t}^{m_0}\notin G(J^{m_0t})$ for
some positive integer $m_0$. From which it  clearly follows that
$f_{j_1}^{l_{1}}\cdots f_{j_t}^{l_{t}}\notin G(J^{m})$ for all
$l_{i}\geq m_0$ where  $m=\sum l_{i}$ .

Notice that in order to show that $f^m \notin G(J^{mt})$ for some $m$, it is enough to show that
$f^k \notin G(\overline{J^{kt}})$ for some $k$. As let $f^k \notin G(\overline{J^{kt}})$ for some $k$. Then,
$f^k=gh$ where $h \in G(\overline{J^{kt}})$ and $\deg g > 0$. Now as $h \in G(\overline{J^{kt}})$, $h^{k_1} \in J^{ktk_1}$
for some $k_1$ which implies $f^{kk_1}=g^{k_1}h^{k_1} \notin G(J^{ktk_1})$.
Hence taking $m=kk_1$, we have $f^m \notin G(J^{mt})$.

We have assumed that $Z\not\subset F$ for any compact face $F \in
\mathcal{F}$, but nevertheless $Z$ may be a subset of a noncompact
face in $\mathcal{F}$. We divide the proof in two cases depending
on whether $Z$ is a subset of some noncompact face or not.

Case 1: First we assume that $Z \not\subset F$
for any face (compact or noncompact) $F \in \mathcal{F}$. Suppose
$f^m \in G(\overline{J^{mt}})$ for all $m$. Without loss of
generality, let $x_1|f$. Since $f \in G(\overline{J^t})$, $g={f}/{x_1} \notin \overline{J^t}$.
Hence $f\in \conv(J^t)$
and $g
\notin \conv(J^t)$. Let $l$ be the line segment joining
$\Gamma(f)$ and $\Gamma(g)$. Then $l$ intersects $\conv({J^t})$ at
some point $p \in tF$ where $F$ is a face of $\conv(J)$, see
Corollary \ref{multiplefaces}. Notice that $p \neq \Gamma(f)$, see Remark
\ref{extra}. Hence, $\Gamma(f)=p+v$ where $0<\|v\|<1$. Now for any
$m$, consider the line segment joining $\Gamma(f^m)$ and
$\Gamma(g^m)$, we denote this line segment by $ml$. We have
$\Gamma(f^m)=mp+mv$ where $mp \in mtF$ and $mtF$ is a face of
$\conv(J^{mt})$. Again as $f^m \in G(\overline{J^{mt}})$,
${f^m}/{x_1} \notin \overline{J^{mt}}$. Notice that $\Gamma(
{f^m}/{x_1})$ and $mp$ lie on $ml$, and since $\Gamma({f^m}/{x_1})
\notin \conv(J^{mt})$ and $mp \in \conv(J^{mt})$, we have
$\|mv\|=\|mp-\Gamma(f^{m})\| \leq
\|\Gamma(f^m)-\Gamma(f^m/x_1)\|=1$ for any $m$, a contradiction.

Case 2: Now assume that $Z \subset G$ for some
noncompact face $G \in \mathcal{F}$ and that
$\{a_{j_1},\ldots,a_{j_t}\} \not\subset F$ for any compact face
$F\in \mathcal{F}$. We prove that $f^m \notin
G(\overline{J^{mt}})$ for some $m=m_0$ by induction on $\dim
G$. If $\dim G=1$, then $f \notin G(\overline{J^{t}})$, because it
follows from Lemma \ref{imp} that the only point on $tG$ which
corresponds to a generator of $\overline{J^{t}}$, is an extremal
point of $\conv(J^t)$ and certainly $a=a_{j_1}+\cdots+a_{j_t}$ is
not an extremal point of $\conv(J^t)$, see Corollary \ref{last}. Now
let $\dim G=p>1$. We may assume that $\{a_{j_1},\ldots,a_{j_t}\}
\not\subset G'$ for any  proper face $G'$ of $G$. As if
$\{a_{j_1},\ldots,a_{j_t}\} \subset G'$ for some proper face $G'$
of $G$, then $G'$ is a noncompact face of $G$ with $\dim G' < \dim
G$ and we are through by induction.

Let $S:=\{ v \in \RR^n \;|\; \langle v,u \rangle = c\}$ be the
supporting hyperplane of $\conv(J)$ such that $S \cap \conv(J)=G$. Since $G$ is a noncompact face,
there exists $j$ such that
$u(j)=0$. Consider
$a_{\lambda}:=a_{j_1}+\cdots+a_{j_t}-\lambda(0,\ldots,1,\ldots,0)$,
$1$ being at $j$th place, $\lambda \geq 0$. Notice that there
exists $\lambda_0>0$ such that $a_{\lambda_0} \notin \conv(J)$.
Let $l_0$ be the line segment joining $a$ and $a_{\lambda_0}$. As
$a \in l_0 \cap tG$, the intersection of  $l_0$ with $tG$ is a
nonempty  convex set. Let $l=l_0 \cap tG$ be the line segment
joining $a$ and $a_{\lambda'}$ where $a_{\lambda'}$ lies on some
proper face $tG'$ of $tG$ and $\lambda' >0$, as $\dim G' < \dim
G$. Also $a_{\lambda'}< a$, so we have $a= a_{\lambda'}+ w$, with
$ \|w\|=\lambda' >0$. For  any positive integer $m$,
$ma_{\lambda'} \in mtG'$ and
$\|ma-ma_{\lambda'}\|=m\|a-a_{\lambda'}\|=m\|w\|>0$. Let  for
$m=m_0$, $m\|w\| \geq 1$. Then for $m=m_0$, $ma$ and
$ma-(0,\ldots,1,\ldots,0)$ lies on $mtG$, $1$ being at $j$ th
place, so that $\Gamma(f^m/x_j) \in mtG$ which implies
$ f^m/x_j \in \overline{J^{mt}}$ and hence $f^m \not\in
G(\overline{J^{mt}})$ for $m=m_0$.
\end{proof}
\noindent

Let $S=K[x_1,\ldots,x_n,y_1,\ldots,y_r]$ be a bigraded polynomial
ring with $\deg x_i=(1,0)$ and $\deg y_j=(d_j,1)$. Recall $J=(f_1,\ldots,f_r)$ where
$f_j=x^{a_j}$ and $\deg f_j=d_j$. Let $ \phi$
be the surjective homomorphism from $S$ to
$\mathcal{R}(J)=K[x_1,\ldots,x_n,f_1t,\ldots,f_rt]$, given by
$x_i \mapsto x_i$ and $y_j \mapsto f_jt$ so that $S/L \iso
\mathcal{R}(J)$ where the ideal $L$ is  generated by binomials of
the type $g_{1}h_{1}-g_{2}h_{2}$ where $g_{1},g_{2}$ are
monomials in $x_i$ and $h_{1},h_{2}$ are monomials in $y_j$.
Notice that $\deg h_1=\deg h_2$.

Now consider the fiber ring
$\mathcal{F}(J)=\mathcal{R}(J)/\mm\mathcal{R}(J)$ of the ideal $J$
where $\mm=(x_1,\ldots,x_n) \subset A$. Then $\mathcal{F}(J) \iso
S/(L,\mm) \iso T/D$ and hence $\mathcal{F}(J)_{\red} \iso T/{\Rad D}$ where $D$ is the image of the ideal $L$ in
$T=S/\mm$, and $T=K[y_1,\ldots,y_r]$. Let $\psi=\phi\tensor
S/\mm\: T\to {\mathcal F}(J)$ be the induced epimorphism. We have
$D=\Ker \psi$. Notice that the ideal $D$ is generated by monomials
and  homogeneous binomials  in the $y_j$. In fact, if
$g_{1}h_{1}-g_{2}h_{2}$ is a generator of $L$, then its image in
$T$ is a monomial, if one of the $g_i$ belongs to $\mm$, and
otherwise it is a homogeneous binomial. We have the following lemma:

\begin{Lemma} \label{common} Let $b=b_1-b_2 \in D$ be  a homogeneous binomial generator of $D$
with $b_{1}=y_{i_1}^{l_{1}}\cdots y_{i_u}^{l_{u}}$,
$b_{2}=y_{j_1}^{m_{1}}\cdots y_{j_v}^{m_{v}}$ and
$\sum_{i=1}^{u}l_{i}=\sum_{j=1}^{v}m_{j}=t$. If the set
 $\{a_{i_1},\ldots,a_{i_u}\} \subset G$ for some $G \in \mathcal{F}_{c}$, then also the set
 $\{a_{j_1},\ldots,a_{j_v}\} \subset G$ .
\end{Lemma}
\begin{proof}
As $b \in D$, we have $\psi(b)=0$, i.e.\
$\psi(b_{1})=\psi(b_{2})$. Therefore we have
$x^{l_{1}a_{i_1}}\cdots x^{l_{u}a_{i_u}}=x^{m_{1}a_{j_1}}\cdots
x^{m_{v}a_{j_v}}$, and so $\sum_{p=1}^u l_{p}a_{i_p}= \sum_{k=1}^v
m_{k}a_{j_k}$. Let the set $\{a_{i_1},\ldots,a_{i_u}\} \subset
G$ for some $G \in \mathcal{F}_{c}$. We show that
$\{a_{j_1},\ldots,a_{j_v}\} \subset G$. Let $S:= \{ v \in \RR^n
\;|\; \langle v,u \rangle = c\}$, be the supporting hyperplane of
$\conv(J)$ such that $S \cap \conv(J)=G$.

We have $\langle \sum_{k=1}^v m_{k}a_{j_k},u \rangle =\langle \sum_{p=1}^u l_{p}a_{i_p},u \rangle=tc$. Suppose
$\{a_{j_1},\ldots,a_{j_v}\} \not\subset G$, then there exists at
least one $k_0 \in\{1,\ldots,v \}$ such that $a_{j_{k_0}} \notin
G$. Since $\langle a_{j_k},u \rangle\geq c$ for all $k$,  it follows that
$\langle a_{j_{k_0}},u \rangle > c$ which in turn implies that $\langle \sum_{k=1}^vl_{k}a_{j_k},u \rangle>tc$,
a contradiction.
\end{proof}
We denote by ${\mathcal F}_{c}$  the set of compact faces, and by
${\mathcal F}_{mc}$ the set of maximal compact faces of $\conv(J)$.
Let $F\in {\mathcal F}_{mc}$; we set $P_F=(y_j: a_j \notin F)$ and
we denote by $B_F$ the kernel of $\theta_F : K[y_j: a_j \in F] \to
K[F]:=K[f_jt: a_j \in F]$ where $\theta_F(y_j)= f_jt$.

With the notation introduced  we have

\begin{Proposition}
\label{difficult} $\Rad D=(\Sect_{F\in {\mathcal
F}_{mc}}P_F,\sum_{F\in {\mathcal F}_{mc}}B_F)=\Sect_{F\in
{\mathcal F}_{mc}}(P_F,B_F).$
\end{Proposition}

\begin{proof} For the proof we proceed in several steps.

\medskip
1.\ Step: Let $f$ be a monomial in $T$. We claim that $f\in \Rad
D\iff f\in \Sect_{F\in{\mathcal F}_{mc}}P_F$.

\medskip
\noindent We may assume that $f$ is squarefree. So let
$f=y_{j_1}\ldots y_{j_k}$ with $j_1<j_2<\cdots < j_k$ and assume
that $f\in \Rad D$. Then $f^{n_0}\in D$ for some integer $n_0$,
and hence $\psi(f^{n_0})=0.$ This implies that
$x^{n_0a_{j1}}\cdots x^{n_0a_{jk}} \in \mm J^{n_0k}$. Hence
$x^{na_{j1}}\cdots x^{na_{jk}}$ is not a minimal generator of
$J^{nk}$ for any $n \geq n_0$. Now  Theorem \ref{2.4} implies that
$\{ a_{j1},\ldots,a_{jk}\} \not\subset F$ for any compact face $F
\in \mathcal{F}$. This shows that  $f\in \Sect_{F \in
\mathcal{F}_{mc}}P_F$.

Conversely, assume that $f\in \bigcap_{F \in
\mathcal{F}_{mc}}P_F$. Then $\{ a_{j_1},\ldots,a_{j_k}\}
\not\subset F$ for any $F \in \mathcal{F}_{mc}$. This implies that
$\{ a_{j_1},\ldots,a_{j_k}\} \not\subset F$ for any compact face.
From Theorem \ref{2.4} we conclude  that there exists an integer
$m$ such that $(x^{a_{j_1}}\cdots x^{a_{j_k}})^m\in \mm J^{km}$.
Since $\psi(f^m)=(x^{a_{j_1}}\cdots x^{a_{j_k}})^m$ it follows
that $f^m\in D$, and hence $f\in\Rad D$.

\medskip
2.\ Step: $D\subset (\Sect_{F\in{\mathcal F}_{mc}}P_F,
\sum_{F\in{\mathcal F}_{mc}}B_F).$

\medskip
\noindent It follows from the first step that all monomial
generators  in   $D$ belong to the ideal $(\Sect_{F\in{\mathcal
F}_{mc}}P_F, \sum_{F\in{\mathcal F}_{mc}}B_F).$  Now let
$b=b_1-b_2$ be one of the homogeneous binomial generators of $D$
with $b_{1}=y_{i_1}^{l_{1}}\cdots y_{i_u}^{l_{u}}$,
$b_{2}=y_{j_1}^{m_{1}}\cdots y_{j_v}^{m_{v}}$ and
$\sum_{i=1}^{u}l_{i}=\sum_{j=1}^{v}m_{j}=t$. As $b \in D$, we
have $\psi(b)=0$, i.e.\ $\psi(b_{1})=\psi(b_{2})$. Therefore we
have $x^{l_{1}a_{i_1}}\cdots
x^{l_{u}a_{i_u}}=x^{m_{1}a_{j_1}}\cdots x^{m_{v}a_{j_v}}$, and so
$\sum_{p=1}^u l_{p}a_{i_p}= \sum_{k=1}^v m_{k}a_{j_k}$. We show
that $b \in \sum_{F\in{\mathcal F}_{mc}}B_F$, if $b \notin
\Sect_{F\in{\mathcal F}_{mc}}P_F$. In fact, if $b \notin
\Sect_{F\in{\mathcal F}_{mc}}P_F$, then  one of the $b_i$, say
$b_1 \notin \Sect_{F\in{\mathcal F}_{mc}}P_F$. This implies that
$\{a_{11},\ldots,a_{1u}\} \in G$ for some compact face $G \in
\mathcal{F}_{mc}$ and  then from Lemma \ref{common},
$\{a_{21},\ldots,a_{2v}\} \in G$. Hence, $b=b_1-b_2 \in B_G$.

\medskip
3.\ Step: $\sum_{F\in{\mathcal F}_{mc}}B_F \subset D.$
\medskip

\noindent
Notice that $B_F = \Ker\theta_F$ and $D= \Ker\psi$. Certainly, for each $F \in \mathcal{F}_{mc}$,
$\Ker\theta_F \subset \Ker\psi$ and hence
$\sum_{F \in \mathcal{F}_{mc}}B_F \subset D$.

\medskip
4.\ Step: $\bigcap_{F \in \mathcal{F}_{mc}}(P_F,B_F)=(\bigcap_{F
\in \mathcal{F}_{mc}} P_F, \sum_{F \in \mathcal{F}_{mc}}B_F).$
\medskip

\noindent For each $F \in \mathcal{F}_{mc}$, let $Q_F =(P_F,B_F)$,
and let $M=\bigcap_{F \in \mathcal{F}_{mc}} P_F$ and $B=\sum_{F
\in \mathcal{F}_{mc}}B_F$. In order to show  that $(M,B)=\bigcap_{F \in \mathcal{F}_{mc}}Q_F$, we proceed in the following
steps:

(i)\;
First we show $(M,B)\subset \bigcap_{F \in \mathcal{F}_{mc}}Q_F.$
Clearly, for each
$F \in \mathcal{F}_{mc}$, $M \subset Q_F$. Now we also prove
that $B \subset Q_F$ for all $F \in \mathcal{F}_{mc}$. Take $b=b_1-b_2 \in B$
with $b_{1}=y_{i_1}^{l_{1}}\cdots y_{i_u}^{l_{u}}$,
$b_{2}=y_{j_1}^{m_{1}}\cdots y_{j_v}^{m_{v}}$ and
$\sum_{i=1}^{u}l_{i}=\sum_{j=1}^{v}m_{j}=t$. Suppose
that $b \notin B_G$, then we prove  $b \in P_G$. As $b \notin
B_G$, it implies that for one of the $b_i$, say for $b_1$, there
exists $y_{i_p}|b_1$ such that $a_{i_p} \notin G$. Once we show that there
exists also some $k \in \{1,\ldots,v\}$ such that $y_{j_k}|b_2$ and
$a_{j_k} \notin G$, then it will imply that $b_1,b_2 \in P_G$ and hence $b \in P_G$.
Suppose this is not the case, then
$\{a_{j_1},\ldots,a_{j_v}\} \in G $. But then from Lemma $\ref{common}$,  we have
$\{a_{i_1},\ldots,a_{i_v}\} \in G $ which is a contradiction.
Hence we have$(M,B) \subset \bigcap_{F \in\mathcal{F}_{mc}}Q_F$.

(ii) Notice that for each $F \in \mathcal{F}_{mc}$, $Q_F$ is a
prime ideal. Indeed, $Q_F$ being the kernel of the surjective map
$\pi_F: K[y_1,\ldots,y_r] \to K[f_it: a_i\in F]$ given by
$\pi_F(y_j)=f_jt$, if $a_j \in F$ and $\pi_F(y_j)=0$, if $a_j
\notin F$, the assertion follows.

(iii)  We claim that $\{Q_F :F \in \mathcal{F}_{mc}\}$ is the set
of all the minimal prime ideals containing $(M,B)$. Let $P$ be any
prime ideal containing $(M,B)$, then it implies that $P \supset
M=\bigcap_{\mathcal{F}_{mc}}P_F$ and so $P \supset P_G$ for some
$G \in \mathcal{F}_{mc}$. Also, $P \supset B= \sum B_F$. Hence $P
\supset Q_G$.

(iv)  We claim $(M,B)$ is a radical ideal, that is,
$\Rad(M,B)=(M,B).$ This  amounts to prove that for all $Q_F$,
$(M,B)T_{Q_F}=Q_FT_{Q_F}$. Fix $G \in \mathcal{F}_{mc}$, the set
$\{y_i : a_i\in G\} \subset T \backslash Q_G$ and hence all $y_i$
such that $a_i \in G$ are invertible in $T_{Q_G}$. For all $P_F$,
$F \neq G$, there exists at least one $y_j\in P_F$ such that $y_j
\in G$, as otherwise $P_F \subset P_G $ which implies $F \supset
G$, a contradiction. Hence for all $F \neq G$,
$P_FT_{Q_G}=T_{Q_G}$. Therefore we have
$(M,B)T_{Q_G}=(\bigcap_{F\in \mathcal{F}_{mc}} P_F, \sum_{F \in
\mathcal{F}_{mc}}B_F)T_{Q_G} =(P_G,\sum_{F \in
\mathcal{F}_{mc}}B_F)T_{Q_G}=(P_G, B_G)T_{Q_G}=Q_GT_{Q_G}.$

\medskip
\noindent Since by (iii) we have  $\Rad(M,B)= \bigcap_{F \in
\mathcal{F}_{mc}}Q_F$ it follows then that  $(M,B)=\bigcap_{F \in
\mathcal{F}_{mc}}Q_F$. Now by Step 1, Step 2 and  Step 3, one has
$$D \subset (\Sect_{F\in{\mathcal F}_{mc}}P_F,
\sum_{F\in{\mathcal F}_{mc}}B_F) \subset \Rad{D}.$$ Finally by
Step 4, we have $(\Sect_{F\in{\mathcal F}_{mc}}P_F,
\sum_{F\in{\mathcal F}_{mc}}B_F)=\bigcap_{F \in
\mathcal{F}_{mc}}(P_F,B_F)$ which is a radical ideal. Hence we
have $\Rad D=(\Sect_{F\in {\mathcal F}_{mc}}P_F,\sum_{F\in
{\mathcal F}_{mc}}B_F)=\Sect_{F\in {\mathcal F}_{mc}}(P_F,B_F).$
\end{proof}

We denote by $\Min(R)$ the set of minimal prime ideals of a ring
$R$.

\begin{Corollary}\label{observation}
Let $I\subset A$ be a monomial ideal. Then there  is an injective
map
\[
{\mathcal F}_{mc}\to \Min({\mathcal F}(I)).
\]
This map is bijective if $I$ is an extremal ideal.
\end{Corollary}
\begin{proof} Let $J$ be the minimal monomial reduction ideal of
$I$.  Then $J$ is an extremal ideal. From above proposition we
have $\mathcal{F}(J)_{\red} \iso T/\Sect_{F\in {\mathcal
F}_{mc}}(P_F,B_F)$ where $(P_F,B_F)$ is a prime ideal for each $F
\in \mathcal{F}_{mc}$. Hence there is a bijective map $$\rho_1\:
{\mathcal F}_{mc}\to \Min({\mathcal F}(J))$$ given by $F \mapsto
(P_F,B_F)/D$.

As $\mathcal{F}(I)$ is integral over $\mathcal{F}(J)$, for each $P
\in \Min(\mathcal{F}(J))$ there exists a minimal prime ideal $Q
\in \Min(\mathcal{F}(I))$ such that $P=Q \cap \mathcal{F}(J)$.
Therefore there exists an injective map $\rho_2$ from
$\Min(\mathcal{F}(J))$ to $\Min(\mathcal{F}(I))$,  and hence
$\rho=\rho_2 \circ \rho_1\: {\mathcal F}_{mc} \to \Min({\mathcal
F}(I))$ is the desired injective map. Finally, if $I$ is extremal,
then  $I=J$ and $\rho=\rho_1$ is a bijection.
\end{proof}

Next corollary gives us a combinatorial characterization of the
fiber ring of an extremal ideal $J$ to be a domain.

\begin{Corollary}\label{smart}
Let $J=(x^{a_1},\ldots,x^{a_r})$ be an extremal ideal. Then the
following conditions are equivalent:
\begin{enumerate}
\item The fiber ring $\mathcal{F}(J)$ is a domain;
\item The reduced fiber ring $\mathcal{F}(J)_{\red}$ is a domain;
\item $|\mathcal{F}_{mc}|=1$.
\end{enumerate}
\end{Corollary}
\begin{proof}
(1) \implies  (2) is obvious,
and (2) \iff (3) follows from Corollary \ref{observation}.\\
(3) \implies (1): Let $|\mathcal{F}_{mc}|=1$. Then it follows from
Proposition \ref{difficult} that  $\Rad D = (B_F,P_F)$ where $F
\in \mathcal{F}_{mc}$. Notice that as there is only one maximal
compact face $F$, the ideal $P_F$ is the zero ideal. Hence
$(P_F,B_F)= B_F$. Also by Step 3 in the proof of Proposition
\ref{difficult} we have $B_F \subset D $. Therefore we have $\Rad
D= D=B_F$ which is a prime ideal. Hence $\mathcal{F}(J)\iso T/D$
is a domain.
\end{proof}

By  the above corollary  the fiber ring of an extremal ideal $J$
is a domain if and only if there is only one maximal compact faces
of $\conv(J)$. But in general the property of being reduced cannot
be  characterized in terms of combinatorial properties of
$\conv(J)$, as the the following simple example demonstrates:

\begin{Example}{\em  Consider the two  extremal ideals $J_1=(x^6,x^2y,xy^2,y^6)$ and $J_2=(x^8,x^6y,x^2y^7,y^{12})$
in the polynomial ring $A=K[x,y]$. It is easy to see that
$\conv(J_1)$ and $\conv(J_2)$ have the same face lattices.
Nevertheless the fiber ring of the ideal $J_1$ given by
$\mathcal{F}(J_1)\iso K[y_1,y_2,y_3,y_4]/(y_1y_4,y_2y_4,y_1y_3)$
is reduced while the fiber ring of the ideal $J_2$ given by
$\mathcal{F}(J_2) \iso K[y_1,y_2,y_3,y_4]/(y_1y_4, y_2y_4^2,
y_2^2y_4-y_1y_3^2, y_1^2y_3)$ is not reduced.}
\end{Example}

Next we define an inverse system of
semigroup rings  $K[F]$ for $F \in \mathcal{F}_c$(set of compact faces of $\conv(I)$) where
$K[F]=K[f_it\: a_i \in F ]$ with $f_i=x^{a_i}$. For $G
\subset F$, define the ring homomorphism $\pi_{GF}: K[F]
\rightarrow K[G]$, given by $\pi_{GF}(f_it)=f_it$, if $a_i \in G$
and $\pi_{GF}(f_it)=0$, otherwise. Notice that $\pi_{GF}$ is well
defined. To see this, we need to show that if
$f_{i_1}f_{i_2}\cdots f_{i_k}t^k=f_{j_1}f_{j_2}\cdots f_{j_k}t^k$
where
$\{{a_{i_{1}}},\ldots,{a_{i_{k}}}\},\{{a_{j_{1}}},\ldots,{a_{j_{k}}}\}
\subset F $, then $\pi_{GF}(f_{i_1}f_{i_2}\cdots
f_{i_k}t^k)=\pi_{GF}(f_{j_1}f_{j_2}\cdots f_{j_k}t^k)$. If
$\pi_{GF}(f_{i_1}\cdots f_{i_k}t^k)=0$, then  $\{a_{i_1},\ldots,
a_{i_k}\} \not\subset G.$ Since $y_{i_1}\cdots
y_{i_k}-y_{j_1}\cdots y_{j_k}\in D$ it follows from Lemma
\ref{common} that $\{a_{j_1},\ldots, a_{j_k}\} \not\subset G$,
too. Hence $\pi_{GF}(f_{j_1}\cdots f_{j_k}t^k)=0$. On the other hand, if $\pi_{GF}(f_{i_1}\cdots f_{i_k}t^k) \neq 0$,
then $\pi_{GF}(f_{j_1}\cdots f_{j_k}t^k) \neq 0$, and so
$$\pi_{GF}(f_{i_1}\cdots f_{i_k}t^k)=f_{i_1}\cdots
f_{i_k}t^k=f_{j_1}\cdots f_{i_k}t^k=\pi_{GF}(f_{j_1}\cdots
f_{i_k}t^k).$$   Hence $\pi_{GF}(f_{i_1}\cdots
f_{i_k}t^k)=\pi_{GF}(f_{j_1}\cdots f_{j_k}t^k)$ in both cases.

Also we may notice that  for $H \subset G \subset F$ and  $F
\in \mathcal{F}_c$, one has  $\pi_{HG}\circ \pi_{GF}= \pi_{HF}.$
Hence the inverse system is well defined.

\begin{Theorem} $F(J)_{\red}\iso\underleftarrow{\lim}_{F \in \mathcal{F}_c}K[F].$
\label{limit}
\end{Theorem}
\begin{proof} For each $F \in \mathcal{F}_c$
consider the ring homomorphism $\pi_F$ from $K[y_1,\ldots,y_r]$ to
$K[F]$ given by $\pi_F(y_j)=f_jt$, if $a_j \in F$ and
$\pi_F(y_j)=0$, if $a_j \notin F$.

Notice that $\Ker\pi_F$ is equal to the ideal $Q_F:=(B_F, P_F)$.
We define the map $$\pi: K[y_1,\ldots y_r] \longrightarrow
\Dirsum_{F \in \mathcal{F}_c}K[F],$$ given by
$\pi=(\pi_F)_{F\in\mathcal{F}_c}$. We have $\Ker\pi= \bigcap_{F
\in \mathcal{F}_c}Q_F=\bigcap_{F \in \mathcal{F}_c}(B_F,P_F)$. We
claim that for all $G \subset F$ one has  $Q_F \subset Q_G$.
Indeed, for all $G \subset F$, $P_F \subset P_G$ and by the proof
of   Proposition \ref{difficult}, Step 4(i), we have $B_F \subset
(B_G,P_G)$. It follows that
$$\Ker\pi=\bigcap_{F \in \mathcal{F}_{mc}}Q_F.$$
Therefore  Proposition \ref{difficult} implies that $\Ker\pi =\Rad
D$.  Thus we have $$K[y_1,\ldots,y_r]/\Ker\pi \iso F(J)_{\red}.$$

It remains to show that $\Im(\pi)=\underleftarrow{\lim}_{F \in
\mathcal{F}_c}K[F]$. First notice that $\mbox{Im} (\pi) \subset
\underleftarrow{\lim}_{F \in \mathcal{F}_c}K[F]$, since
$\pi_{GF}\circ \pi_F=\pi_G$ for all $G\subset F$.

Now let $v=(m_F)_{F \in \mathcal{F}_c} \in
\underleftarrow{\lim}_{F \in \mathcal{F}_c}K[F]$. We may assume
that for each $F\in \mathcal{F}_c$, the element $m_F$ is a
monomial in $K[F]$ since all homomorphisms in the inverse system
are multigraded. For each $F \in \mathcal{F}_c$, we choose $g_F
\in K[y_1,\ldots,y_r]$ such that $\pi_F(g_F)=m_F$ and with the
property that whenever $m_F=m_G$ in $K[x_1,\ldots,x_n,t]$ then it
implies $g_F=g_G$. (Notice that for each $F \in \mathcal{F}$, the
$K$-algebra  K[F] can be naturally embedded in the $K$-algebra
$K[x_1,\ldots,x_n,t]$).

Let $Z=\{m_F :m_F \neq 0, F \in
\mathcal{F}_c\}=\{m_1,\ldots,m_l\}$. For each $i=1,\ldots ,l$, we
define the set $A_i=\{F \in \mathcal{F}_c: m_F=m_i \}$. We claim
that for each $A_i$ one has $\bigcap_{F \in A_i} F \in A_i$. Fix
an $i$, and notice that it is enough to show that for any $F,G \in
A_i$ we have $F \cap G \in A_i$. Let $m_F=f_{i_1}\cdots
f_{i_p}t^p=f_{j_1}\cdots f_{j_p}t^p=m_G$. Then it follows by Lemma
\ref{common} that the sets $\{a_{i_1},\ldots,a_{i_p}\},
\{a_{j_1},\ldots,a_{j_p}\} \subset F \cap G =H$. Therefore
$\pi_{HF}(m_F)=m_F$ and $\pi_{HG}({m_G})=m_G$. Also as $v=(m_F)_{F
\in \mathcal{F}_c} \in \underleftarrow{\lim}_{F \in
\mathcal{F}_c}K[F]$ we have $\pi_{HF}(m_F)=m_H=\pi_{HG}(m_G)$.
Hence $m_G=m_F=m_H$, so $H \in A_i$. Hence $H_i= \bigcap_{F \in
A_i} F \in A_i$, $i=1,\ldots,l$.

For each $i$, we choose a monomial $g_{H_i} \in K[y_1,\ldots,y_r]$
such that $\pi_{H_i}(g_{H_i})=m_{H_i}$. For all $F \in A_i$, we
define $g_F=g_{H_i}$, $i=1,\ldots,l$ and for all $F \in
\mathcal{F}_c \backslash \bigcup_{i=1}^lA_i$, we define $g_F=0$.
Notice that for all $F \in \mathcal{F}_c$, we have
$\pi_F(g_F)=m_F$. Indeed, let $F \in \mathcal{F}_c$. If $F \in
\mathcal{F}_c \backslash \bigcup_{i=1}^lA_i$, then $g_F=0=m_F$ and
we have $\pi_F(g_F)=m_F$. If $F \in A_i$ for some $i$, then as we
have $\pi_{H_iF}\circ \pi_F=\pi_{H_i}$ and
$\pi_{H_i}(g_F)=m_{H_i}=m_F$, it follows by the very definition of
the map $\pi_{H_iF}$ that $\pi_F(g_F)=m_F$. Moreover, by our
choice of the $g_F$ we also have $g_F=g_G$ whenever $m_F=m_G$.

Now let $S=\{g_F : F\in\mathcal{F}_{mc}\}$, and let $g=\sum_{g_F
\in S}g_F$. We claim that $\pi(g)=v$, i.e. $\pi_G(g)=m_G$ for all
$G \in \mathcal{F}_{c}$. Notice that it is enough to show that
$\pi_G(g)=m_G$ for all $G \in \mathcal{F}_{mc}$.  In fact, if $H
\in \mathcal{F}_{c}$ there exists  $G \in \mathcal{F}_{mc}$ such
that $H \subset G$, and  since $\pi_{G}(g)=m_G$, we have
$\pi_{H}(g)=\pi_{HG}(\pi_{G}(g))=\pi_{HG}(m_G)=m_H$.

Now let $G\in \mathcal{F}_{mc}$. We claim that $\pi_G(g_F)=0$ for
all $g_F\neq g_G$, so that we have $\pi_G(g)=m_G$, as asserted.

To prove this claim, let $g_F=y_{i_1}\cdots y_{i_p}$ and suppose
that $\pi_G(g_F) \neq 0$. Then we have
$\{a_{i_1},\ldots,a_{i_p}\} \subset G \cap F$. Let $H=G \cap F$,
then $H \in \mathcal{F}_c$. Since $v \in \underleftarrow{\lim}_{F
\in \mathcal{F}_c}K[F]$ and $H$ is a common face of $F$ and $G$,
we have $\pi_{HF}(m_F)=m_H=\pi_{HG}(m_G)$. As
$\{a_{i_1},\ldots,a_{i_p}\} \subset H$, we have $0 \neq
m_F=\pi_{HF}(m_F)=m_H=\pi_{HG}(m_G)=m_G$.
Hence $g_F=g_G$, a contradiction.
\end{proof}

The analytic spread $\ell$ of any ideal $I$ in a Noetherian local
ring $(R,\mm)$ is given by the Krull dimension of the fiber ring
$\mathcal{F}(I)$ of $I$. It has been shown by Carles Bivia-Ausina
\cite{BA} that the analytic spread of any non-degenerate ideal $I
\subset \CC[[x_1,\ldots,x_n]]$ is equal to the $c(I)+1$ where
$$c(I)=\max\{\dim{F} : F\; \mbox{is a compact face of} \conv(I)\}.$$

Next we show that for monomial ideals this result is an immediate
consequence of our structure theorem (Theorem \ref{limit}).
\begin{Corollary}
\label{ausina} Let $I \subset A=K[x_1,\ldots,x_n]$ be any monomial
ideal. Let $\ell = \dim \mathcal{F}(I)$ be the analytic spread of
ideal $I$. Then
$$\ell= c(I)+1=\max\{\dim F: F\;  \mbox {is a compact face of}\; \conv(I)\}+1.$$
\end{Corollary}
\begin{proof}
Let $J$ be the minimal monomial reduction ideal of $I$.  We have
$\ell=\dim
\mathcal{F}(I)=\dim\mathcal{F}(J)=\dim\mathcal{F}(J)_{\red}$. By
Theorem \ref{limit}, we have $\mathcal{F}(J)_{\red}
=\underleftarrow{\lim}_{F \in \mathcal{F}_c}K[F]\subset \Dirsum_{F
\in \mathcal{F}_c}K[F]$. Therefore $\dim(\mathcal{F}(J))\leq
\max\{\dim K[F]: F \in \mathcal{F}_c\}.$ As $\dim K[F]=\dim F+1$,
it follows that  $\ell \leq c(I)+1$.

\noindent For proving $\ell \geq c(I)+1$, we notice that the
canonical homomorphisms $$\bar{\pi}_G\: \underleftarrow{\lim}_{F
\in \mathcal{F}_c}K[F]\to K[G]$$ are surjective for all $G\in
\mathcal{F}_{c}$.
\noindent Indeed, if $m$ is a monomial in $K[G]$
and $v=(m_F)_{F\in \mathcal{F}_{c}} \in \underleftarrow{\lim}_{F
\in \mathcal{F}_c}K[F]$ with
\[
m_F
=\left\{ \begin{array}{lll} m, & \mbox{if} & \supp(m)\subset F,\\
0, & \mbox{if} & \supp(m)\subset F,
\end{array} \right.
\]
then  $\bar{\pi}_F(v)=m$. Here $\supp(m)$ of some monomial
$m=x_1^{a_1}\cdots x_n^{a_n}\in A$  is defined to be
$\supp(m)=\{a_i : a_i \neq 0\}$.

\noindent It follows that $ \dim F(J) \geq \dim K[F]$ for all $F
\in \mathcal{F}_{mc}$. Therefore we have $\ell \geq c(I)+1$, as
desired.
\end{proof}

\section{On the reduction number of a monomial ideal}

In this section we consider the reduction number of a monomial
ideal $I \in A$ with respect to the minimal monomial reduction
ideal $J$. We show in Corollary \ref{normal} that if $I^m$ is
integrally closed for $m \leq \ell$ then $I$ is normal and the
reduction number of $I$ with respect to $J$ is less than $\ell-1$.
Here $\ell$ denotes the analytic spread of the monomial ideal $I$
and the reduction number of an ideal $I$ with respect to $J$ is
defined to be the minimum of $m$ such that $JI^m=I^{m+1}$.

\begin{Theorem} \label{ell} Let $I \subset A=K[x_1,\ldots,x_n]$ be a monomial ideal and $J$
its minimal monomial reduction ideal. Let $\ell$
be the analytic spread of $I$. Then
$$\overline{I^{m}}=J\overline{I^{m-1}}\quad \text{for all}\quad m\geq
\ell.$$
\end{Theorem}
\begin{proof} We may assume $I$ is a proper ideal, and let $I=(x^{a_1},x^{a_2},\ldots,x^{a_s})$  where
$f_i=x^{a_i}=x_1^{a_i(1)}x_2^{a_i(2)}\cdots x_n^{a_i(n)}$ for
$i=1,\ldots,s$. Without loss of generality, let $J= (
x^{a_1},x^{a_2},\ldots,x^{a_r})$ be the minimal monomial reduction
ideal of $I$ so that $\ext(I)=\{a_1,\ldots,a_r\}$. Let $m \geq
\ell$,
 we show $\overline{I^m} \subset J \overline{I^{m-1}}$, the other inclusion being trivial.
Let $x^b\in \overline{I^m}=\overline{J^m}$ where $x^b=x_1^{b(1)}\cdots x_n^{b(n)}$.

For the proof we consider the following two cases:
\medskip

Case 1. $b \in F$ where $F$ is a face of $\conv(I^m)$.
\medskip

First we claim that $b=b_1+v$ where $b_1 \in G $ for some compact
face $G$ of $\conv(I^m)$ and $v \in \RR^{n}_{\geq 0}$. If $F$ is a
compact face, then we take $v=0$ and $b_1=b$. Now let $F$ be  a
noncompact face. We prove the claim by induction on $\dim F$. If
$\dim F =1$, then clearly $b=ma_i+ v$ where $v \in \RR^{n}_{\geq
0}$ for some $a_i \in \ext(I)$. Now let $\dim F=t>1$. Let $S =\{ v
\in \RR^n \;|\; \langle v,u \rangle = c \}$ (where $
u=(u(1),\ldots,u(n))\in \RR^n$, $c\in \RR$)  be a supporting
hyperplane of $\conv(I^m)$ such that $S \cap \conv(I^m)=F.$ Since
$F$ is an noncompact face there exists $u(j)$ such that $u(j)=0$.
Consider $b_{\lambda}:=b-\lambda(0,\ldots,1,\ldots,0)$, $1$ being
at $j$th place, $\lambda \geq 0$. Notice that there exists
$\lambda_0>0$ such that $b_{\lambda_0} \notin \conv(I^m)$. Let
$l_0$ be the line segment joining $b$ and $b_{\lambda_0}$. The
intersection of $l_0$ with $F$, is nonempty and therefore is a
convex set. It follows that $l=l_0 \cap F$ is a line segment
joining $b$ and $b_{\lambda'}$ where $b_{\lambda'}$ lies on some
proper face $F'$ of $F$ and $\lambda'\geq 0$. Therefore $bb_{\lambda'}+ w$ with $b_{\lambda'} \in F'$ and $w\in \RR^n_{\geq
0}$. By induction, $b_{\lambda'}=b_1 + w'$ where $b_1 \in G$ for
some compact face $G$ and $w' \in \RR^n_{\geq 0}$. Hence $b=b_1+v$
with $v=w+w' \in \RR^n_{\geq 0}$. Hence the claim.

As $G$ is a compact face we have $\dim G \leq \ell$ by Corollary
\ref{ausina}. Now since $b_1 \in G$ and there exists $p\leq \ell$
affinely independent vectors $\{a_{i_1},\ldots,a_{i_p}\} \subset
\ext(I)$ such that $b_1=\sum_{j=1}^p k_ja_{i_j}$ with $\sum k_i
=m$. Since $p \leq \ell \leq m$, there exists $a_{i_{j_0}}$ such
that $b_1-a_{i_{j_0}} \in \conv(I^{m-1})$. Therefore,
$b-a_{i_{j_0}}=b_1-a_{i_{j_0}}+v \in \conv(I^{m-1}) \cap
\NN^n=\Gamma(\overline{I^{m-1}})$. Hence $b \in
\Gamma(J\overline{I^{m-1}})$.

\medskip

Case 2. $b \notin F$ for any face $F$ of $\conv(I^m)$.
\medskip

\noindent Let $f=x^b$. We may assume that $f \in
G(\overline{J^m})$. Without loss of generality, let $x_1|f$. Since
$f \in G(\overline{J^m})$, $g= {f}/{x_1} \notin \overline{J^m}$.
Hence $b \in \conv(I^m)$ and $\Gamma(g) \notin \conv(I^m)$. Let
$l$ be the line segment joining $b$ and $\Gamma(g)$. Then $l$
intersects $\conv({I^m})$ at some point $a \in F$ where $F$ is a
face of $\conv(I^m)$. Hence, $b=a+v$ where $ v \in \RR^n_{\geq
0}$. Now by the proof of first case, we may write $a=a_1+v_1$ such
that $a_1 \in G$ for some compact face $G$ of $\conv(I^m)$ and $
v_1 \in \RR^n_{\geq 0}$. Hence $b=a_1+w$ where $w=v+v_1 \in
\RR^n_{\geq 0}$. Hence as in the above case, we get that $x^b \in
J \overline{I^{m-1}}$.
\end{proof}

\Remark{\em There is a related result by Wiebe. He shows that for
the maximal graded ideal $\mm$ in a positive normal affine
semigroup ring $S$ of dimension $d$ one has
$\overline{\mm^{n+1}}=\mm\overline{\mm^{n}}$ for all $n \geq d-2$,
and that $\overline{\aa^{n+1}}=\aa\overline{\aa^{n}}$ for all $n
\geq d-1$ if  $\aa \subset S$ is an integrally closed ideal, see
\cite[Theorem 2.1]{A}. }

\begin{Corollary} \label{normal} Let $I^{a}$ be integrally closed for all $a \leq \ell-1$, then
$I^\ell=JI^{\ell-1}$ and $I$ is normal, i.e. $I^{a}$ is integrally closed for all $a$.
\end{Corollary}

\begin{proof} By the above theorem we have $\overline{I^{\ell}} \subset J \overline{I^{\ell-1}}$, and since $\overline{I^{\ell-1}}=I^{\ell-1}$,
we see that $\overline{I^{\ell}} \subset JI^{\ell-1}$. Hence
$I^\ell=JI^{\ell-1}$.

\noindent
Also, $\overline{I^{\ell}}= J \overline{I^{\ell-1}}=JI^{\ell-1} \subset I^{\ell} \subset \overline{I^{\ell}}$.
Hence $\overline{I^\ell}=I^{\ell}$. By applying induction on $a$, one has $\overline{I^a}=I^a$ for all $a$.
\end{proof}

\Remarks{\em

(a) Corollary \ref{normal}\; is a generalization of a result by
Reid, Roberts and Vitulli \cite[Proposition 2.3]{RRV} . They
proved that if $I \subset A=K[x_1,\ldots,x_n]$ is a monomial ideal
and $I^m$ is integrally closed for $m \leq n-1$, then $I$ is a
normal ideal.

(b) In Corollary \ref{normal}, once we assume that the monomial
ideal $I$ is normal, then the bound on the reduction number with
respect to monomial reductions can be obtained as a consequence of
a theorem by  Valabrega-Valla  \cite{VV} and the improved version
of the Briancon-Skoda theorem due to Aberbach and Huneke
\cite{AH}. Infact, if $I$ is a normal monomial ideal, then $R(I)$
is Cohen-Macaulay and hence $F(I)$ is Cohen-Macaulay. Thus by
Valabrega-Valla  \cite{VV} and Aberbach-Huneke \cite{AH}, the
reduction number of $I$ with respect to monomial reductions is
less than the analytic spread $\ell$ of $I$. I am thankful to
Prof.\ Verma for this information.
 }

\end{document}